\documentclass[nopreprintline,sort&compress]{elsarticle}

\usepackage{lipsum}
\makeatletter
\def\ps@pprintTitle{%
 \let\@oddhead\@empty
 \let\@evenhead\@empty
 \def\@oddfoot{}%
 \let\@evenfoot\@oddfoot}
\makeatother

\usepackage{amsmath,amssymb}
\usepackage{multirow}
\usepackage{graphicx}
\usepackage[hidelinks]{hyperref}
\usepackage{caption}
\usepackage{subcaption}

\newcommand{\p}{\partial}
\newcommand{\erf}{\operatorname{erf}}

\begin{document}

\title{Numerical simulations with the finite element method for the Burgers' equation on the real line}

\author[dmpa]{P.H.A. Konzen\corref{cor}}
\ead{pedro.konzen@ufrgs.br}

\author[dmpa,ppgmap]{E. Sauter}
\ead{esequia@gmail.com}

\author[dmpa,ppgmap]{F.S. Azevedo}
\ead{fazedo@gmail.com}

\author[dmpa,ppgmap]{P. Zingano}
\ead{pzingano@mat.ufrgs.br}

\address[dmpa]{Departamento de Matem\'{a}tica Pura e Aplicada,
     IME - Instituto de Matem\'{a}tica e Estat\'{i}stica,
     UFRGS - Universidade Federal do Rio Grande do Sul, 91509-900 Porto Alegre, RS, Brasil}

\address[ppgmap]{PPGMAp - Programa de P\'{o}s-Gradua\c{c}\~{a}o em Matem\'{a}tica Aplicada}

\begin{abstract}
In this paper we present a simple and accurate second order finite element scheme to simulate the Burgers' equation on the whole real line and subjected to initial conditions with compact support. The numerical simulations are performed by considering a sequence of auxiliary spatially dimensionless Dirichlet's problems parameterized by the domain's semidiameter $L$. Gaining advantage from the well-known convective-diffusive effects of the Burgers' equation, computations start by choosing $L$ larger than the semidiameter of the support of the initial condition and, as solution diffuses out, $L$ is increased appropriately. By direct comparisons between numerical and analytic solutions and its asymptotic behavior, we conclude this simple scheme is very accurate and can be applied to numerically investigate properties of this and similar equations on infinite domains.
\end{abstract}


\maketitle

\section{Introduction}

Consider the viscous Burger's equation defined on the real line:
\begin{equation}
  \label{eq:burgers}
  \frac{\p u}{\p t} + u\frac{\p u}{\p x} = \nu \frac{\p^2 u}{\p x^2},\quad (x\in\mathbb{R}, t>0),
\end{equation}
subjected to the initial condition:
\begin{equation}
  \label{eq:initial_condition}
  u(x,0) = g(x),\quad (x\in\mathbb{R}),
\end{equation}
where $\nu > 0$ is a given viscosity coefficient and $g$ is a given function with compact support on $\mathbb{R}$.

Burgers' equation is known to have appeared firstly in $1915$ in the work of Harry Bateman \cite{Bateman1915a}, but it receives its name after to the Dutch physicist J.M. Burgers, who applied this equation in the understanding of turbulent fluids \cite{Burgers1974a}. This homogeneous quasilinear parabolic partial differential equation appears in the modeling of several phenomena such as shock flows, wave propagation in combustion chambers, vehicular traffic movement, acoustic transmission, etc. (see, for instance, \cite{Fletcher1982a} and the references therein). Another import characteristic of this equation is its several well known analytic solutions in bounded and unbounded domains. Therefore, this equation is already a classical test case in mathematical analysis and numerical simulations of convective-diffusive partial differential equations.

From the analytic point of view the literature is rich in discussing solutions and properties for the Burgers' equation on bounded and unbounded regions and subjected to a variety of initial and boundary conditions (see, for instance, ~\cite{Bastos2006a,Burgers1974a, Abd-el-Malek2000a, Evans2010a, Gorguis2006a, Holland1977a, Rodin1970a, Wood2006a}). Now, from the numerical simulation point of view the majority of the studies found in the literature are concerned about the Burgers' equation defined in a bounded region and subjected to Dirichlet's boundary conditions. Several numerical schemes have been applied to simulate this problem, for instance: Finite Element Methods \cite{Aksan2005a,Aksan2006a,Arminjon1981a,Caldwell1982a,Caldwell1981a,Caldwell1987a,Dogan2004a,Fletcher1983a,Hrymak1986a,Kadalbajoo2006a,Kutluay2004a,Ozis2003a,Ozis2005b,Shao2011a,Zhang2009a}, Finite Difference Methods \cite{Caldwell1982a,Fletcher1983a,Gulsu2006a,Kutluay1999a,Mukundan2015a}, variational schemes \cite{Aksan2004a,Caldwell1985a,Ozis1996a}, spectral methods \cite{Basdevant1986a,Khater2008a}, Hardy's multiquadric method \cite{Hon1998a}, matched asymptotic expansion methods \cite{Ozis2005a}, multisymplectic box methods \cite{Tabatabaei2007a}, Homotopy Analysis Methods \cite{Inc2008a}, the quintic B-spline collocation procedure \cite{Saka2008a}, the gradient reproducing kernel particle method \cite{Hashemian2008a}, quasi-interpolation techniques \cite{Xu2011a}, uniform Haar wavelets \cite{Jiwari2015a}.

In this work we present an efficient numerical scheme based on the Finite Element Method to simulate Burgers' equation on the real line and subjected to initial conditions with compact support. The proposed scheme explore the convective-diffusive nature of the differential equation. If for small times the convective effects are predominant demanding very fine and localized meshes, for large times diffusion takes place and the solution tend to relax demanding less refined but large meshes. We deal with it by computing the finite element discretization of the dimensionless spatially form of the Burgers' equation on a fixed mesh and, then, allowing the diameter of the domain to vary as the solution spreads out. This simple idea has been proved very computational efficient producing accurate results. This is supported by direct comparisons between numerical and analytic solutions and their asymptotic behavior.

In the next section we briefly discuss the analytic solution of problem \eqref{eq:burgers}-\eqref{eq:initial_condition} and its asymptotic properties. In Section \ref{sec:finite_element_scheme} we present the proposed time and space discretization of the spatially dimensionless form of the Burgers' equation. In Section \ref{sec:implementation_scheme} we discuss the details of the implementation scheme. Then in Section \ref{sec:numerical_experiments} we present numerical experiments, which endorse the efficiency and accuracy of the scheme as to its potential to be applied to investigate solution properties on the real line. Finally, in Section \ref{sec:final_considerations} we close by summarizing the principal aspects of this work.

\section{Analytic solution}\label{sec:analytic_solution}

Here we recall the well known closed-form expression for $u(x,t)$ obtained by J. Cole and E. Hopf \cite{Cole1951a, Hopf1950a}. Introducing $\theta(x,t)$ and $\theta_{0}(x)$ by the Hopf-Cole transformation:
\begin{equation}
  \theta(x,t) := \exp\left\{-\frac{1}{2\nu}\int_0^x u(y,t)\,dy\right\},\quad \theta_0(x) := \exp\left\{-\frac{1}{2\nu}\int_0^x g(y)\,dy\right\}
\end{equation}
one obtains that $\theta$ solves the following initial value problem for the heat equation:
\begin{align}
  &\frac{\p \theta}{\p t}  = \nu \frac{\p^{2} \theta}{\p x^{2}},\quad (x \in \mathbb{R}, t>0)\\
  &\theta(x,0) = \theta_0(x),\quad (x \in \mathbb{R}),
\end{align}
whose unique bounded solution is given by Poisson's formula:
\begin{equation}
\theta(x,t) = \frac{1}{\sqrt{4 \pi\nu t}}\int_{-\infty}^\infty e^{-\frac{|x-y|^2}{4\nu t}}\theta_0(y)\,dy,\quad (x\in\mathbb{R}, t>0).
\end{equation}
Since $\displaystyle u = -2 \nu \frac{\theta_{x}}{\theta}$, it follows that:
\begin{equation}\label{eq:analytic_solution}
  u(x,t) = \frac{\int_{-\infty}^\infty \frac{x-y}{t}e^{-\frac{|x-y|^2}{4\nu t}}\theta_0(y)\,dy}{\int_{-\infty}^\infty e^{-\frac{|x-y|^2}{4\nu t}}\theta_0(y)\,dy},\quad (x\in\mathbb{R}, t>0).
\end{equation}

This also shows that problem \eqref{eq:burgers}-\eqref{eq:initial_condition} has a unique solution $\displaystyle u(\cdot,t) \in C^{0}\left([0, \infty)\right.$, $\left. L^{1}(\mathbb{R})\right)$, given by \eqref{eq:analytic_solution} above, which satisfies: $u \in C^{\infty}(\mathbb{R} \times (0, \infty))$ and $u(\cdot,t) \in C^{0}\left((0,\infty), W^{k,p}(\mathbb{R})\right)$ for every $k \geq 1$, $p \geq 1$. Here, $W^{k,p}(\mathbb{R})$ is the Sobolev space of functions in $L^{p}(\mathbb{R})$ whose $k$-th order derivatives belong to $ L^{p}(\mathbb{R}) $. Moreover, by \eqref{eq:analytic_solution} and standard heat kernel estimates one gets that: 
\begin{align}
  &\|u(\cdot,t)\|_{L^{2}(\mathbb{R})} = O(t^{-\frac{1}{4}}), &&\|u(\cdot,t)\|_{L^{\infty}(\mathbb{R})} = O(t^{-\frac{1}{2}}), \\
  &\|u_{x}(\cdot,t)\|_{L^{2}(\mathbb{R})} = O(t^{-\frac{3}{4}}), &&\|u_{x}(\cdot,t)\|_{L^{\infty}(\mathbb{R})} = O(t^{-1}), \\
  &\|u_{xx}(\cdot,t)\|_{L^{2}(\mathbb{R})} = O(t^{-\frac{5}{4}}), &&\|u_{xx}(\cdot,t)\|_{L^{\infty}(\mathbb{R})} = O(t^{-\frac{3}{2}}),
\end{align}
and so on.

A more refined analysis in \cite{Zingano1997a} shows that the asymptotic limits:
\begin{equation}\label{eq:asymptote}
\gamma_{p} \equiv \lim_{t \to \infty} t^{\frac{1}{2}\left(1 - \frac{1}{p}\right)}\|u(\cdot,t)\|_{L^p(\mathbb{R})}, \quad 1 \leq p \leq \infty,
\end{equation}
are well defined and have the following values. Let $m$ be the solution mass, that is:
\begin{equation}
m = \int_{-\infty}^\infty u(x,t)\,dx = \int_{-\infty}^\infty u_0(x)\,dx.
\end{equation}
Then:
\begin{equation}\label{eq:gamma_p_values}
\gamma_{p} = \frac{|m|}{\sqrt{4\pi\nu}} (4\nu)^{\frac{1}{2 p}}\frac{2 \nu}{m}\left(1 - e^{-\frac{m}{2\nu}}\right) \|\mathcal{F}\|_{L^p(\mathbb{R})}
\end{equation}
with ${\cal F} \in L^1(\mathbb{R})\cap L^\infty(\mathbb{R})$ defined by:
\begin{equation}
  \mathcal{F}(x) = \frac{e^{-x^2}}{\lambda - h \erf(x)}
\end{equation}
where $\erf(\cdot)$ is the error function:
\begin{equation}
  \erf(x) = \frac{2}{\sqrt{\pi}}\int_0^x e^{-\xi^2}\,d\xi
\end{equation}
and $\lambda$, $h$ are given by:
\begin{equation}
  \lambda = \frac{1 + e^{-\frac{m}{2\nu}}}{2},\quad h = \frac{1 - e^{-\frac{m}{2\nu}}}{2}.
\end{equation}

When $p = 1$, \eqref{eq:gamma_p_values} is simply:
\begin{equation}
\lim_{t\to \infty} \|u(\cdot,t)\|_{L^1(\mathbb{R})} = |m|,
\end{equation}
and we further have: $\displaystyle \|u_{x}(\cdot,t)\|_{L^{1}(\mathbb{R})} = O(t^{-\frac{1}{2}})$, $\displaystyle \|u_{xx}(\cdot,t)\|_{L^{1}(\mathbb{R})} = O(t^{-1})$, and so on.

These results will be used in Section  \ref{sec:numerical_experiments} as further evidence for the accuracy of the numerical approximation scheme developed in the next two sections.

\section{Finite element scheme}\label{sec:finite_element_scheme}

We consider the following auxiliary Dirichlet's problem:
\begin{align}
  & \frac{\p \tilde{u}}{\p t} + \frac{1}{L}\tilde{u}\frac{\p \tilde{u}}{\p \tilde{x}} = \frac{\nu}{L^2}\frac{\p^2 \tilde{u}}{\p \tilde{x}^2},\quad (x\in (-1, 1), t>0), \label{eq:burgers_dimensionless}\\
  & \tilde{u}(\tilde{x},0) = \tilde{g}(\tilde{x}),\quad (x\in (-1, 1)),\label{eq:burgers_dimensionless_ic}\\
  & \tilde{u}(-1,t) = \tilde{u}(1,t) = 0,\quad (t>0), \label{eq:burgers_dimensionless_bc}
\end{align}
where $\tilde{x} := x/L$ is the dimensionless space variable, $L$ is the reference semidiameter of the domain and $\tilde{g}(\tilde{x}) := g(L\cdot\tilde{x})$. From now one we will work with this space dimensionless problem and, for the sake of simplicity, we will omit the tilde, i.e., we will denotes $\tilde{x}$ simply by $x$ and $\tilde{u}$ by $u$.

Following the Rothe's method, we start by discretizing equation \eqref{eq:burgers_dimensionless} in time. To this end, we consider the following $\theta$-scheme for the time discretization of equation \eqref{eq:burgers_dimensionless}:
\begin{equation}\label{eq:burgers_discrete_in_time}
  \begin{split}
  \frac{u^{n+1} - u^{n}}{\Delta t} &=  - \frac{\theta}{L} u^{n+1}\frac{\p u^{n+1}}{\p x} - \frac{(1 - \theta)}{L} u^{n}\frac{\p u^{n}}{\p x} \\
  &+ \frac{\theta}{L^2} \nu \frac{\p^2 u^{n+1}}{\p x^2} + \frac{(1 - \theta)}{L^2} \nu \frac{\p^2 u^{n}}{\p x^2}
  \end{split}
\end{equation}
where $u^0 = u(x,0)$, $u^{n}$ denotes the approximation of $u(x,t_n)$, $n=1, 2, \dotsc$, $t_n = n\Delta t$, $\Delta t$ is a given time step size and $0 \leq \theta \leq 1$. For simplicity sake, from now one we denote $u^{n+1}$ by $u$ and $u^{n}$ by $u^0$.

Now, we consider the following weak formulation of the problem defined by equations \eqref{eq:burgers_discrete_in_time}, \eqref{eq:burgers_dimensionless_ic} and \eqref{eq:burgers_dimensionless_bc}: given $u^0\in H^1_0(-1, 1)$ find $u\in H^1_0(-1, 1)$ such that:
\begin{equation}\label{eq:weak_formulation}
  \begin{split}
    &\left(\varphi, u\right) + \frac{\theta}{L} \Delta t \left(\varphi, u\frac{\p u}{\p x}\right) + \frac{\theta}{L^2} \Delta t \left(\frac{\p \varphi}{\p x}, \nu \frac{\p u}{\p x}\right) - \\
    &\left(\varphi, u^0\right) + \frac{(1 - \theta)}{L} \Delta t \left(\varphi, u^0\frac{\p u^0}{\p x}\right) + \frac{(1 - \theta)}{L^2} \Delta t \left(\frac{\p \varphi}{\p x}, \nu \frac{\p u^0}{\p x}\right) = 0    
  \end{split}
\end{equation}
for all $\varphi\in H_0^1(-L, L)$.

Let's consider the following second order finite element triple $(\mathcal{K}, P_2(\mathcal{K}), \Sigma)$, where the cells $\mathcal{K}\subset\mathcal{T}_h$ are line segments forming a regular triangulation $\mathcal{T}_h$ of the segment $[-1, 1]$, the element shape functions $P_2(\mathcal{K}) = \{v:\mathcal{K}\to \mathbb{R}, v(x) = a_0 + a_1x + a_2x^2, a_0,a_1,a_2\in\mathbb{R}\}$ are second order polynomials, and the degrees of freedom $\Sigma$ are located at the end points of each $\mathcal{K}$ and its middle point (see, for instance, \cite{Johnson2009a}). This allows us to define the finite element space:
\[
V_h := \{v\in C^0(-1, 1): v|_{\mathcal{K}}\in P_2(\mathcal{K}), \forall \mathcal{K}\in\mathcal{T}_h\} \subset H^1_0[-1, 1].
\]

Then, following the Galerkin's method, we iteratively approximate the solution of \eqref{eq:burgers_dimensionless} subjected to \eqref{eq:burgers_dimensionless_ic} and \eqref{eq:burgers_dimensionless_bc} by the solution of the following full discrete problem: given $u_h^0\in V_h$ find $u_h\in V_h$ such that:
\begin{equation}
  \label{eq:discrete_problem}
  \begin{split}
    &\left(\varphi_i, u_h\right) + \frac{\theta}{L} \Delta t \left(\varphi_i, u_h\frac{\p u_h}{\p x}\right) + \frac{\theta}{L^2} \Delta t \left(\frac{\p \varphi_i}{\p x}, \nu \frac{\p u_h}{\p x}\right) - \\
    &\left(\varphi_i, u_h^0\right) + \frac{(1 - \theta)}{L} \Delta t \left(\varphi_i, u_h^0\frac{\p u_h^0}{\p x}\right) + \frac{(1 - \theta)}{L^2} \Delta t \left(\frac{\p \varphi_i}{\p x}, \nu \frac{\p u_h^0}{\p x}\right) = 0
  \end{split}
\end{equation}
for all $\varphi_i$ in the basis of the finite element space $V_h$.

At each time step, we solve the nonlinear system of equations \eqref{eq:discrete_problem} by the Newton's method. The Newton's formulation than reads: given $u_h^{0}\in V_h$ we iteratively compute approximations $u_h^{m+1}$ of $u_h$ by iterating:
\begin{subequations}\label{eq:newton_iteration}
  \begin{equation}
    \label{eq:newton_system}
    J(u_h^{(m)})\delta u^{(m)} = -F(u_h^{(m)}) 
  \end{equation}
  \begin{equation}
    u_h^{(m+1)} = u_h^{(m)} + \delta u^{(m)}
  \end{equation}
\end{subequations}
where $F(u_h^{m})$ denotes the left-hand-side of equation \eqref{eq:discrete_problem} substituting there $u_h$ by $u_h^{m}$, $\delta u^{m}$ is the Newton update, and the Jacobian matrix $J(u) = [\mathfrak{j}_{i,j}]_{i,j=0}^{M,M}$ have its elements defined by:
\begin{equation}
  \label{eq:jacobian_elements}
  \begin{split}
    \mathfrak{j}_{i,j} &:= \left(\varphi_i, \varphi_j\right) + \frac{\theta}{L} \Delta t \left(\varphi_i, \frac{\p u}{\p x}\varphi_j\right) + \frac{\theta}{L} \Delta t \left(\varphi_i, u \frac{\p \varphi_j}{\p x}\right) \\
    &+ \frac{\theta}{L^2} \Delta t \left(\frac{\p \varphi_i}{\p x}, \nu \frac{\p \varphi_j}{\p x}\right)    
  \end{split}
\end{equation}
where $M+1$ counts for the number of degrees of freedom.

We now lead the discussion to the implementation of this standard finite element scheme to simulate the Burgers' equation defined on the real line and subjected to an initial condition with compact support. 

\section{Implementation scheme} \label{sec:implementation_scheme}

Because of the convective-diffusive nature of the Burgers' equation, very fine meshes are demanded to accurately compute the solution for small times, but as time increases the solution tend to relax allowing the application of less refined meshes. By assuming an initial condition with compact support numerical simulations of the auxiliary Dirichlet's problem \eqref{eq:burgers_dimensionless}-\eqref{eq:burgers_dimensionless_bc} may produce accurate solutions for finite times. To ensure the accuracy we just need to choose appropriate time and spatial meshes, and pick $L$ sufficiently large. However, the larger the physical time we would like to consider the larger $L$ should be.

The convective effects are predominant for small times and it is appropriate to work with a small $L$, which reduces the demanding in the number of vertices of the discretization scheme. On the other hand, as time increases, the solution spreads out demanding a larger $L$, but a less refined mesh, to ensure the accuracy of the numerical simulation. We deal with this paradigm as follows.

Let's $K$ be the compact support of the initial condition, i.e. $K := \{x\in\mathbb{R}:~g(x)\neq 0\}$ is a compact subset of $\mathbb{R}$. Without loss of generality, we assume that $0\in K$ and denote $d = \max_{x\in K} \{|x|\}$. Also, let's denote by $\mathcal{K}_\Gamma$ the boundary elements of the finite element space $(\mathcal{K}, P_2(\mathcal{K}), \Sigma)$. With this in mind, the implementation idea is to start simulating the auxiliary Dirichlet's problem \eqref{eq:burgers_dimensionless}-\eqref{eq:burgers_dimensionless_bc} by choosing an appropriate $L > d$. Then, at each time iteration $n$ we check if the numerical support $\tilde{K} := \{x\in\Sigma:~|u_h^{n}(x)| > 10^{-15}\}$ is still a subset of $\mathcal{T}_h\setminus \{\mathcal{K}_\Gamma\}$. If it is not the case, then we simply increase $L$ and interpolate the numerical solution $u_h^{n}$ onto $\mathcal{T}_{2h}$.

We point out that the above implementation scheme does not demand one to rewrite the finite element triangulation at each increasing of the reference parameter $L$, since we are always simulating using the same fixed triangulation built in the domain $[-1, 1]$. Moreover, if we exactly double $L$ when $\tilde{K}$ is no longer a subset of $\mathcal{T}_h\setminus \{\mathcal{K}_\Gamma\}$, then the interpolation that is done to restart the computations with the new $L$ is performed by a simple relocation of elements of the coordinate solution vector.

We summarize the implementation procedure as follows:
\begin{enumerate}
\item Set a uniform mesh with $N$ vertices built in the domain $[-1, 1]$.
\item Set the finite element triangulation $\mathcal{T}_h$.
\item Set an appropriate $L > d$.
\item Set the initial solution vector $u_h^0 \leftarrow [g(x_i)]_{i=0}^{2N}$, where $x_i$ is the abscissa of the $i$-th degree of freedom.
\item Set the present solution vector $u_h \leftarrow u_h^0$.
\item Set the time step $\delta_t$.
\item Loop over time steps:
  \begin{enumerate}
  \item Loop over Newton steps:
    \begin{enumerate}
    \item Assemble the Newton system.
    \item Solve the system.
    \item Set $u_h \leftarrow u_h + \delta u_h$.
    \end{enumerate}
  \item If the numerical support of the computed solution is not a subset of $\mathcal{T}_h\setminus \{\mathcal{K}_\Gamma\}$, then:
    \begin{enumerate}
    \item Set $L \leftarrow 2L$.
    \item Rearrange the elements of vector $u_h$.
    \end{enumerate}
  \item Set $u^0 \leftarrow u$.
  \end{enumerate}
\end{enumerate}

The numerical simulations where implemented in C++ using the deal.II open source finite element library \citep{dealII83}. We applied the UMFPACK sparse direct linear solver implemented there to compute the Newton update $\delta u^m$ from equation \eqref{eq:newton_system}. Evaluations of the analytic solution and its asymptotic behavior were performed in Python using the numerical quadrature available in the Scipy module for integration \cite{Scipy}.

\section{Numerical experiments} \label{sec:numerical_experiments}

Here, we present numerical simulations of problem \eqref{eq:burgers} subjected to the initial condition \eqref{eq:initial_condition} with:
\begin{equation}
  \label{eq:numerical_experiment_ic}
  g(x) = \left\{\begin{array}{ll}e^{-10x^2} &, -2 \leq x \leq 2\\
    0 &, \text{otherwise}\end{array} \right.
\end{equation}
and several values of the diffusion coefficient $\nu$.

We first present direct comparisons between the numerical and analytic solution. Then, in the Subsection \ref{subsec:solution_properties} we show that the proposed numerical scheme is able to preserve important analytic properties of the Burgers' equation.

\subsection{Numerical {\it versus} analytic solutions}

 The analytic solution \eqref{eq:analytic_solution} was evaluated with a precision of five significant digits. All the reported simulations were performed applying the Cranck-Nicolson method by choosing the parameter $\theta = \frac{1}{2}$ in the finite element scheme. The stop criteria for the Newton iterations was chosen to be $\|\delta u^{(m)}\| < 10^{-10}$, where $\|\cdot\|$ denotes the $l_2$-vector norm.

Table \ref{tab:nu=1.0-tf=0.05} presents numerical and analytic results of the problem \eqref{eq:burgers}-\eqref{eq:initial_condition} with $\nu = 1$ at time $t_f = 0.05$ and at the domain points $x = -1.0, -0.5, 0.0, 0.5, 1.0$. Numerical results obtained by finite element simulations with different mesh sizes and time steps are reported. We point out that the proposed numerical scheme has provided results in very good agreement with the analytic solution. The relative accuracy of the numerical results is at least $0,02\%$ for this case. Simulations with meshes of $N = 401$ and $N = 801$ vertices have been sufficient to produce a precision of $5$ significant digits. Moreover, the simulations with time step $\Delta t = 10^{-4}$ provided a small gain in accuracy against the simulations with $\Delta t = 10^{-3}$. These good characteristics of the obtained numerical approximations can be also observed in Table \ref{tab:nu=1.0-tf=0.5-100.0}, which reports results at times $t_f = 0.5, 2.5, 10.0, 100.0$ and at several different domain points. 

\begin{figure}
    \centering
    \begin{subfigure}[b]{0.45\textwidth}
    \centering
        \includegraphics[width=\textwidth]{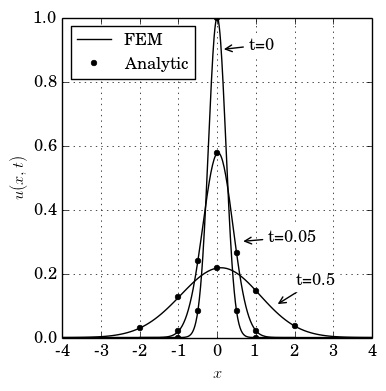}
        \caption{}
    \end{subfigure}
    ~
    \begin{subfigure}[b]{0.45\textwidth}
    \centering
        \includegraphics[width=\textwidth]{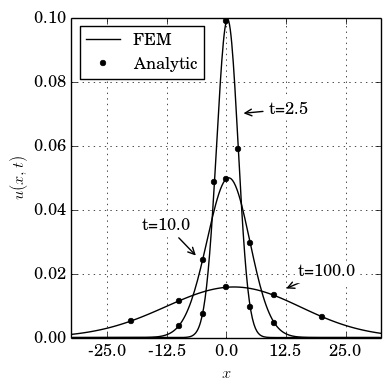}
        \caption{}
    \end{subfigure}
    \caption{Numerical {\it versus} analytic results for the Burgers' equation with $\nu = 1.0$. (a) $t = 0.0, 0.05, 0.5$. (b) $t = 2.5, 10.0, 100.0$.}\label{fig:ab-nu=1.0-dt-1e-3-vert=401}
\end{figure}

Graphical comparisons between numerical and analytic solutions when $\nu = 1.0$ are found in Figure \ref{fig:ab-nu=1.0-dt-1e-3-vert=401} . With this large diffusion coefficient we can observe a small effect given by convection, while solution diffuses rapidly and its numerical support increases in diameter. 

\begin{table}
  \centering
  \caption{Solutions for $\nu = 1.0$ at $t_f = 0.05$.}
  \label{tab:nu=1.0-tf=0.05}
  \begin{tabular}{rrccccc} \hline
    $\Delta t$ & $x$    & $N=101$ & $N=201$ & $N=401$ & $N=801$ & Analytic \\ \hline
    \multirow{5}{*}{$10^{-3}$}   & $-1.0$ & $1.9932$ & $1.9933$ & $1.9933$  & $1.9933$ & $1.9935(-02)$ \\
               & $-0.5$ & $2.3848$    & $2.3850$ & $2.3850$ & $2.3850$ & $2.3849(-01)$ \\
               &  $0.0$ & $5.7619$    & $5.7619$ & $5.7619$ & $5.7619$ & $5.7621(-01)$ \\
               &  $0.5$ & $2.6430$    & $2.6432$ & $2.6432$ & $2.6432$ & $2.6432(-01)$ \\
               &  $1.0$ & $2.1310$    & $2.1311$ & $2.1311$ & $2.1311$ & $2.1314(-02)$ \\ \hline
    \multirow{5}{*}{$10^{-4}$}   & $-1.0$ & $1.9934$ & $1.9935$ & $1.9935$ & $1.9935$    \\
               & $-0.5$ & $2.3848$ & $2.3849$ & $2.3849$ & $2.3849$   \\
               &  $0.0$ & $5.7620$ & $5.7620$ & $5.7621$ & $5.7621$  \\
               &  $0.5$ & $2.6429$ & $2.6432$ & $2.6432$ & $2.6432$  \\
               &  $1.0$ & $2.1312$ & $2.1314$ & $2.1314$ & $2.1314$  \\ \hline
  \end{tabular}
\end{table}

\begin{table}
  \centering
  \caption{Solutions for $\nu = 1.0$ at $t_f = 0.5, 2.5, 10.0, 100.0$, FEM simulations with $\Delta t = 10^{-3}$.}
  \label{tab:nu=1.0-tf=0.5-100.0}
  \begin{tabular}{rrccccc} \hline
    $\Delta t$ & $x$    & $N=101$ & $N=201$ & $N=401$ & $N=801$ & Analytic \\ \hline
    \multicolumn{7}{c}{$t_f = 0.5$} \\ \hline
    \multirow{5}{*}{$10^{-3}$}   & $-2.0$ & $2.9490$ & $2.9476$ & $2.9476$    & $2.9476$ & $2.9476(-02)$ \\
               & $-1.0$ & $1.2535$    & $1.2538$ & $1.2539$ & $1.2539$ & $1.2539(-01)$ \\
               &  $0.0$ & $2.1720$    & $2.1720$ & $2.1720$ & $2.1720$ & $2.1720(-01)$ \\
               &  $1.0$ & $1.4615$    & $1.4620$ & $1.4621$ & $1.4621$ & $1.4621(-01)$ \\
               &  $2.0$ & $3.5973$    & $3.5960$ & $3.5960$ & $3.5960$ & $3.5960(-02)$ \\ \hline
    \multirow{5}{*}{$10^{-4}$}   & $-2.0$ & $2.9490$ & $2.9476$ & $2.9476$ & $2.9476$ \\
               & $-1.0$ & $1.2535$ & $1.2538$ & $1.2539$  & $1.2539$\\
               &  $0.0$ & $2.1720$ & $2.1720$ & $2.1720$  & $2.1720$ \\
               &  $1.0$ & $1.4615$ & $1.4620$ & $1.4621$  & $1.4621$\\
               &  $2.0$ & $3.5973$ & $3.5960$ & $3.5960$  & $3.5960$\\ \hline
    \multicolumn{7}{c}{$t_f = 2.5$} \\ \hline
    \multirow{5}{*}{$10^{-3}$}   & $-5.0$ & $7.4475$   & $7.4533$ & $7.4539$ & $7.4538$ & $7.4538(-03)$ \\
               & $-2.5$ & $4.8764$    & $4.8752$ & $4.8751$ & $4.8750$ & $4.8750(-02)$ \\
               &  $0.0$ & $9.8941$    & $9.8942$ & $9.8942$ & $9.8942$ & $9.8942(-02)$ \\
               &  $2.5$ & $5.8835$    & $5.8818$ & $5.8816$ & $5.8815$ & $5.8815(-02)$ \\
               &  $5.0$ & $9.4493$    & $9.4557$ & $9.4564$ & $9.4563$ & $9.4563(-03)$ \\ \hline
    \multirow{5}{*}{$10^{-4}$}   & $-5.0$ & $7.4475$ & $7.4533$ & $7.4539$ & $7.4538$ \\
               & $-2.5$ & $4.8764$ & $4.8752$ & $4.8751$ & $4.8750$ \\
               &  $0.0$ & $9.8941$ & $9.8942$ & $9.8942$ & $9.8942$ \\
               &  $2.5$ & $5.8835$ & $5.8818$ & $5.8816$ & $5.8815$ \\
               &  $5.0$ & $9.4493$ & $9.4557$ & $9.4564$ & $9.4563$ \\ \hline
    \multicolumn{7}{c}{$t_f = 10.0$} \\ \hline
    \multirow{5}{*}{$10^{-3}$}   & $-10.0$ & $3.6372$ & $3.6401$  & $3.6404$ & $3.6404$ & $3.6404(-03)$ \\
               &  $-5.0$ & $2.4244$   & $2.4239$ & $2.4238$ & $2.4237$ & $2.4237(-02)$ \\
               &   $0.0$ & $4.9634$   & $4.9635$ & $4.9635$ & $4.9635$ & $4.9635(-02)$ \\
               &   $5.0$ & $2.9520$   & $2.9512$ & $2.9511$ & $2.9510$ & $2.9510(-02)$ \\
               &  $10.0$ & $4.6962$   & $4.6994$ & $4.6997$ & $4.6997$ & $4.6997(-03)$ \\ \hline
    \multirow{5}{*}{$10^{-4}$}   & $-10.0$ & $3.6372$ & $3.6402$ & $3.6404$ & $3.6404$ \\
               &  $-5.0$ & $2.4244$ & $2.4239$ & $2.4238$ & $2.4237$ \\
               &   $0.0$ & $4.9634$ & $4.9635$ & $4.9635$ & $4.9635$ \\
               &   $5.0$ & $2.9520$ & $2.9512$ & $2.9511$ & $2.9510$ \\
               &  $10.0$ & $4.6962$ & $4.6994$ & $4.6998$ & $4.6997$ \\ \hline
    \multicolumn{7}{c}{$t_f = 100.0$} \\ \hline
    \multirow{5}{*}{$10^{-3}$}   & $-20.0$ & $5.1825$ & $5.1822$ & $5.1822$ & $5.1822$ & $5.1822(-03)$ \\
    & $-10.0$ & $1.1419$ & $1.1418$ & $1.1418$ & $1.1418$ & $1.1418(-02)$ \\
    &   $0.0$ & $1.5709$ & $1.5709$ & $1.5709$ & $1.5709$ & $1.5709(-02)$ \\
    &  $10.0$ & $1.3179$ & $1.3179$ & $1.3179$ & $1.3179$ & $1.3179(-02)$ \\
    &  $20.0$ & $6.5373$ & $6.5366$ & $6.5366$ & $6.5366$ & $6.5366(-03)$ \\ \hline
  \end{tabular}
\end{table}

Table \ref{tab:dt-1e-3-vert=801-L} contains the physical times when the semidiameter $L$ were doubled in performing numerical simulations with time step $\Delta t = 10^{-3}$ and mesh of $801$ vertices for different diffusion coefficients $\nu$. To enhance the discussion about the accuracy of the results, we will give attention to solution's profiles obtained with different values of $L$.

\begin{table}
  \centering
  \caption{Physical times when the semidiameter $L$ were doubled in performing the numerical simulations with time step $\Delta t = 10^{-3}$ and mesh of $801$ vertices with different diffusion coefficients $\nu$.}
  \label{tab:dt-1e-3-vert=801-L}
  \begin{tabular}{lllll} \hline
    $L$          & $\nu = 1.0$         & $\nu = 0.1$         & $\nu = 0.01$        & $\nu=0.001$\\ \hline
    $2\to 4$     & $5.00\times 10^{-3}$ & $5.00\times 10^{-3}$ & $4.98\times 10^{-1}$ & $4.55\times 10^0$\\
    $4\to 8$     & $9.80\times 10^{-2}$ & $9.76\times 10^{-1}$ & $8.05\times 10^0$   & $1.59\times 10^1$\\
    $8\to 16$    & $4.76\times 10^{-1}$ & $4.70\times 10^{0}$  & $3.06\times 10^{1}$ & $5.83\times 10^1$\\
    $16\to 32$   & $2.02\times 10^0$   & $1.96\times 10^{1}$  & $1.19\times 10^{2}$ & $2.21\times 10^2$\\
    $32\to 64$   & $8.35\times 10^0$   & $8.00\times 10^{1}$  & $4.69\times 10^{2}$ & $$\\ 
    $64\to 128$  & $3.43\times 10^{1}$ & $3.26\times 10^{2}$   & $$                 & $$\\ 
    $128\to 256$ & $1.41\times 10^{2}$ & $$   & $$\\ \hline
  \end{tabular}
\end{table}

\begin{figure}
    \centering
    \begin{subfigure}[b]{0.45\textwidth}
    \centering
        \includegraphics[width=\textwidth]{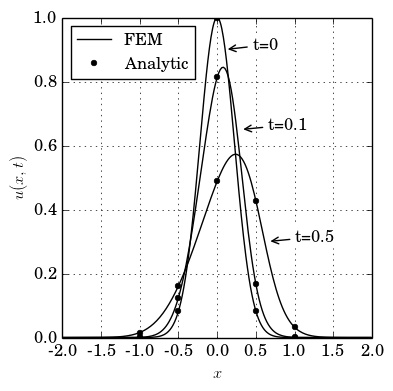}
        \caption{}
    \end{subfigure}
    ~
    \begin{subfigure}[b]{0.45\textwidth}
    \centering
        \includegraphics[width=\textwidth]{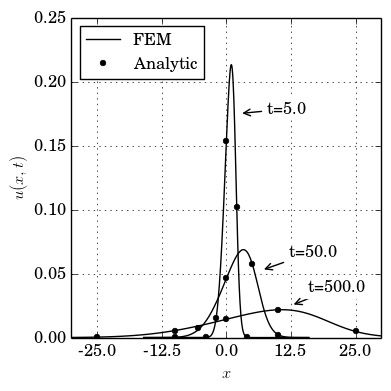}
        \caption{}
    \end{subfigure}
    \caption{Numerical {\it versus} analytic results for the Burgers' equation with $\nu = 0.1$. (a) $t = 0.0, 0.1, 0.5$. (b) $t = 5.0, 50.0, 500.0$.}\label{fig:ab-nu=0.1-dt-1e-3-vert=401}
\end{figure}

Table \ref{tab:nu=0.1-dt-1e-3} reports numerical and analytic results when $\nu = 0.1$ at times $t_f = 0.1, 1.0, 5.0, 50.0, 500.0$ and at several different points $x$. As before, results show very good agreement with at least four significant digits. These results are also graphically presented in Figure \ref{fig:ab-nu=0.1-dt-1e-3-vert=401}. We observe that the convective effects are stronger than before, once we now set a smaller diffusion coefficient.

\begin{table}
  \centering
  \caption{Solutions for $\nu = 0.1$ at $t_f = 0.1, 1.0, 5.0, 50.0, 500.0$. FEM simulations with $\Delta t = 10^{-3}$.}
  \label{tab:nu=0.1-dt-1e-3}
  \begin{tabular}{rccccc} \hline
    $x$    & $N=201$ & $N=401$ & $N=801$ & Analytic \\ \hline
    \multicolumn{5}{c}{$t_f = 0.1$} \\ \hline
    $-1.0$ & $6.6382$ & $6.6379$ & $6.6379$  & $6.6379(-04)$ \\
    $-0.5$ & $1.2484$ & $1.2484$ & $1.2484$  & $1.2484(-01)$ \\
     $0.0$ & $8.1288$ & $8.1289$ & $8.1289$  & $8.1289(-01)$ \\
     $0.5$ & $1.6601$ & $1.6601$ & $1.6601$  & $1.6601(-01)$ \\
     $1.0$ & $6.7263$ & $6.7258$ & $6.7258$  & $6.7257(-04)$ \\ \hline
    \multicolumn{5}{c}{$t_f = 1.0$} \\ \hline
    $-2.0$ & $1.2237$ & $1.2236$ & $1.2236$  & $1.2236(-04)$ \\
    $-1.0$ & $3.6493$ & $3.6493$ & $3.6493$  & $3.6493(-02)$ \\
     $0.0$ & $3.5397$ & $3.5397$ & $3.5397$  & $3.5397(-01)$ \\
     $1.0$ & $1.3624$ & $1.3624$ & $1.3624$  & $1.3624(-01)$ \\
     $2.0$ & $2.1259$ & $2.1256$ & $2.1256$  & $2.1256(-04)$ \\ \hline    
    \multicolumn{5}{c}{$t_f = 5.0$} \\ \hline
    $-4.0$ & $6.0530$ & $6.0527$ & $6.0526$  & $6.0526(-05)$ \\
    $-2.0$ & $1.4916$ & $1.4916$ & $1.4916$  & $1.4916(-02)$ \\
     $0.0$ & $1.5387$ & $1.5387$ & $1.5387$  & $1.5387(-01)$ \\
     $2.0$ & $1.0178$ & $1.0178$ & $1.0178$  & $1.0178(-01)$ \\
     $4.0$ & $3.0283$ & $3.0280$ & $3.0280$  & $3.0280(-04)$ \\ \hline        
    \multicolumn{5}{c}{$t_f = 50.0$} \\ \hline
    $-10.0$ & $1.9051$ & $1.9048$ & $1.9048$  & $1.9048(-04)$ \\
     $-5.0$ & $7.8304$ & $7.8305$ & $7.8305$  & $7.8305(-03)$ \\
      $0.0$ & $4.6189$ & $4.6189$ & $4.6189$  & $4.6189(-02)$ \\
      $5.0$ & $5.7506$ & $5.7505$ & $5.7505$  & $5.7505(-02)$ \\
     $10.0$ & $2.2609$ & $2.2606$ & $2.2606$  & $2.2606(-03)$ \\ \hline
    \multicolumn{5}{c}{$t_f = 500.0$} \\ \hline
    $-25.0$ & $3.4511$ & $3.4513$ & $3.4513$  & $3.4512(-04)$ \\
    $-10.0$ & $5.4509$ & $5.4509$ & $5.4509$  & $5.4509(-03)$ \\
      $0.0$ & $1.4289$ & $1.4289$ & $1.4289$  & $1.4289(-02)$ \\
     $10.0$ & $2.1700$ & $2.1701$ & $2.1701$  & $2.1701(-02)$ \\
     $25.0$ & $4.7812$ & $4.7812$ & $4.7812$  & $4.7812(-03)$ \\ \hline
  \end{tabular}
\end{table} 

\begin{figure}
    \centering
    \begin{subfigure}[b]{0.45\textwidth}
    \centering
        \includegraphics[width=\textwidth]{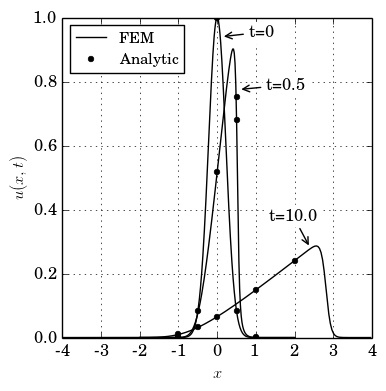}
        \caption{}
    \end{subfigure}
    ~
    \begin{subfigure}[b]{0.45\textwidth}
    \centering
        \includegraphics[width=\textwidth]{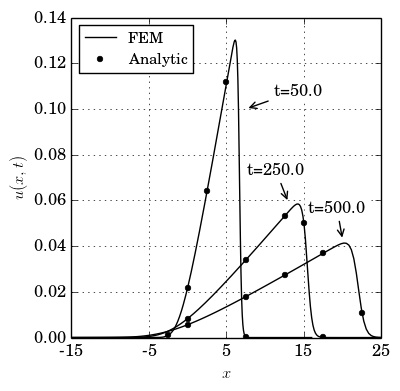}
        \caption{}
    \end{subfigure}
    \caption{Numerical {\it versus} analytic results for the Burgers' equation with $\nu = 0.01$. (a) $t = 0.0, 0.5, 10.0$. (b) $t = 50.0, 250.0, 500.0$.}\label{fig:ab-nu=0.01-dt-1e-3-vert=401}
\end{figure}

Table \ref{tab:nu=0.01-dt-1e-3} presents numerical {\it versus} analytic solutions for the case when $\nu = 0.01$ at times $t_f = 0.5, 10, 50, 250, 500$ and at several different points. This smaller diffusion coefficient $\nu = 0.01$ cases the convective effects to be much more strong even at large times (see, also, Figure \ref{fig:ab-nu=0.01-dt-1e-3-vert=401}). Nevertheless, diffusion still big enough to our numerical scheme to produce very accurate simulations with a mesh of $801$ vertices and time step of $\Delta t = 10^{-3}$.

\begin{table}
  \centering
  \caption{Solutions for $\nu = 0.01$ at $t_f = 0.5, 10.0, 50.0, 250.0, 500.0$, FEM simulations with $\Delta t = 10^{-3}$.}
  \label{tab:nu=0.01-dt-1e-3}
  \begin{tabular}{rccc} \hline
    $x$    & $N=401$ & $N=801$ & Analytic \\ \hline
    \multicolumn{4}{c}{$t_f = 0.5$} \\ \hline
    $-1.0$ & $2.1788$ & $2.1788$  & $2.1788(-04)$ \\
    $-0.5$ & $7.5111$ & $7.5111$  & $7.5111(-02)$ \\
     $0.0$ & $5.1787$ & $5.1787$  & $5.1787(-01)$ \\
     $0.5$ & $6.8116$ & $6.8111$  & $6.8111(-01)$ \\
     $1.0$ & $2.2105$ & $2.2105$  & $2.2105(-04)$ \\ \hline
    \multicolumn{4}{c}{$t_f = 10.0$} \\ \hline
    $-1.0$ & $9.5488$ & $9.5488$  & $9.5488(-03)$ \\
    $-0.5$ & $3.1517$ & $3.1517$  & $3.1517(-02)$ \\
     $0.0$ & $6.5267$ & $6.5267$  & $6.5267(-02)$ \\
     $1.0$ & $1.4914$ & $1.4914$  & $1.4914(-01)$ \\
     $2.0$ & $2.4069$ & $2.4069$  & $2.4069(-01)$ \\ \hline
    \multicolumn{4}{c}{$t_f = 50.0$} \\ \hline
    $-2.5$ & $1.0140$ & $1.0140$  & $1.0140(-03)$ \\
     $0.0$ & $2.1888$ & $2.1888$  & $2.1888(-02)$ \\
     $2.5$ & $6.3993$ & $6.3993$  & $6.3993(-02)$ \\
     $5.0$ & $1.1119$ & $1.1119$  & $1.1119(-01)$ \\
     $7.5$ & $3.2072$ & $3.2128$  & $3.2127(-04)$ \\ \hline
    \multicolumn{4}{c}{$t_f = 250.0$} \\ \hline
     $0.0$ & $8.3037$ & $8.3037$  & $8.3037(-03)$ \\
     $7.5$ & $3.3935$ & $3.3935$  & $3.3935(-02)$ \\
    $12.5$ & $5.3122$ & $5.3122$  & $5.3122(-02)$ \\
    $15.0$ & $4.9921$ & $4.9922$  & $4.9922(-02)$ \\
    $17.5$ & $5.8448$ & $5.8370$  & $5.8379(-05)$ \\ \hline
    \multicolumn{4}{c}{$t_f = 500.0$} \\ \hline
     $0.0$ & $5.6266$ & $5.6266$  & $5.6266(-03)$ \\
     $7.5$ & $1.7910$ & $1.7910$  & $1.7910(-02)$ \\
    $12.5$ & $2.7264$ & $2.7264$  & $2.7264(-02)$ \\
    $17.5$ & $3.6903$ & $3.6903$  & $3.6903(-02)$ \\
    $22.5$ & $1.0871$ & $1.0872$  & $1.0872(-02)$ \\ \hline
  \end{tabular}
\end{table}

\begin{figure}
    \centering
    \begin{subfigure}[b]{0.45\textwidth}
    \centering
        \includegraphics[width=\textwidth]{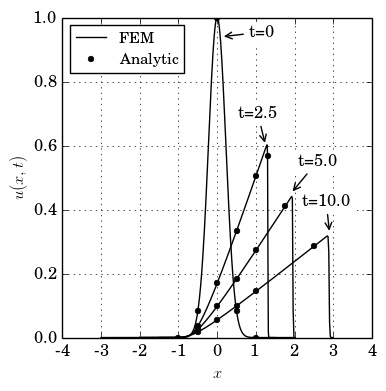}
        \caption{}
    \end{subfigure}
    ~
    \begin{subfigure}[b]{0.45\textwidth}
    \centering
        \includegraphics[width=\textwidth]{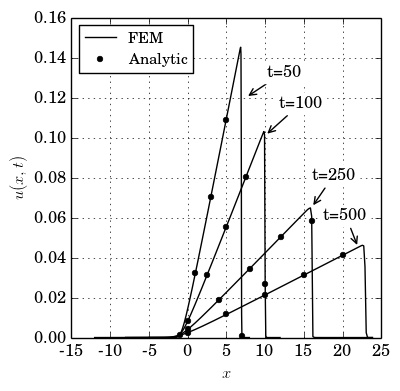}
        \caption{}
    \end{subfigure}
    \caption{Numerical {\it versus} analytic results for the Burgers' equation with $\nu = 0.001$. (a) $t = 0.0, 2.5, 5.0, 10.0$. (b) $t = 50.0, 100.0, 250.0, 500.0$.}\label{fig:ab-nu=0.001-dt-1e-3-vert=401}
\end{figure}

The smaller diffusion coefficient which we will report here is $\nu = 0.001$. Table \ref{tab:nu=0.001-dt-1e-3} presents the numerical {\it versus} analytic solutions for this case at times $t_f = 5, 50, 100, 250$ and at several different points. This small diffusion coefficient provokes solutions with almost a shock wave (see Figure \ref{fig:ab-nu=0.001-dt-1e-3-vert=401}). Even though we can see that our numerical scheme can produce very accurate simulations by taking a time step of $\Delta t = 10^{-3}$ and a mesh with $801$ points. Simulations for smaller diffusion coefficients are also possible, but they demand more (local) refined meshes and smaller time steps.

\begin{table}
  \centering
  \caption{Solutions for $\nu = 0.001$ at $t_f = 5.0, 50.0, 100.0, 250.0$, FEM simulations with $\Delta t = 10^{-3}$.}
  \label{tab:nu=0.001-dt-1e-3}
  \begin{tabular}{rccc} \hline
    $x$    & $N=401$ & $N=801$ & Analytic \\ \hline
    \multicolumn{4}{c}{$t_f = 5.0$} \\ \hline
    $-0.5$ & $2.5377$ & $2.5377$  & $2.5377(-02)$ \\
     $0.0$ & $9.7790$ & $9.7790$  & $9.7790(-02)$ \\
     $0.5$ & $1.8310$ & $1.8310$  & $1.8310(-01)$ \\
     $1.0$ & $2.7253$ & $2.7253$  & $2.7253(-01)$ \\
    $1.75$ & $4.0992$ & $4.0992$  & $4.0992(-01)$ \\ \hline
    \multicolumn{4}{c}{$t_f = 50.0$} \\ \hline
    $-1.0$ & $1.5250$ & $1.5250$  & $1.5250(-03)$ \\
     $1.0$ & $3.2281$ & $3.2281$  & $3.2281(-02)$ \\
     $3.0$ & $7.0537$ & $7.0537$  & $7.0537(-02)$ \\
     $5.0$ & $1.0955$ & $1.0955$  & $1.0955(-01)$ \\
     $7.0$ & $8.3132$ & $6.2548$  & $6.1865(-04)$ \\ \hline
    \multicolumn{4}{c}{$t_f = 100.0$} \\ \hline
     $0.0$ & $8.1649$ & $8.1649$  & $8.1649(-03)$ \\
     $2.5$ & $3.1193$ & $3.1193$  & $3.1193(-02)$ \\
     $5.0$ & $5.5535$ & $5.5535$  & $5.5535(-02)$ \\
     $7.5$ & $8.0129$ & $8.0129$  & $8.0129(-02)$ \\
    $10.0$ & $2.6770$ & $2.6774$  & $2.6781(-02)$ \\ \hline
    \multicolumn{4}{c}{$t_f = 250.0$} \\ \hline
     $0.0$ & $4.0513$ & $4.0513$  & $4.0513(-03)$ \\
     $4.0$ & $1.8749$ & $1.8749$  & $1.8749(-02)$ \\
     $8.0$ & $3.4433$ & $3.4433$  & $3.4433(-02)$ \\
    $12.0$ & $5.0260$ & $5.0260$  & $5.0260(-02)$ \\
    $16.0$ & $6.0637$ & $5.8466$  & $5.8109(-02)$ \\ \hline
  \end{tabular}
\end{table}

We now report error measurements in $L_1$-, $L_2$- and $L_\infty$-norms of the numerical solutions against analytic solutions. Because of the structure of the finite element space and the dimensionless form of the auxiliary problem \eqref{eq:burgers_dimensionless}-\eqref{eq:burgers_dimensionless_ic}, we compute these norms as follows:
\begin{subequations}\label{eq:norms}
  \begin{equation}
    \|e_h\|_{L^1(\mathbb{R})} = L \sum_{\mathcal{K}\in\mathcal{T}_h} \int_\mathcal{K} |e_h(y,t)|\,dy,
  \end{equation}
  \begin{equation}
    \|e_h\|_{L^2(\mathbb{R})} = \left(L \sum_{\mathcal{K}\in\mathcal{T}_h} \int_\mathcal{K} |e_h(y,t)|^2\,dy\right)^{\frac{1}{2}},
  \end{equation}
  \begin{equation}
    \|e_h\|_{L^\infty(\mathbb{R})} = \max_{\mathcal{K}\in\mathcal{T}_h} \{|e_h(y,t)|\},
  \end{equation}
\end{subequations}
where, $e_h(\cdot, t) := u_h(\cdot, t) - u(\cdot, t)$.

\begin{table}[h!]
  \centering
  \caption{Maximum of the error norms of the finite element solution with $\delta_t = 10^{-3}$ against the analytic solution on the time interval $t\in [0, 1]$.}
  \begin{tabular}[h!]{llccc} \hline
    $\nu$ & $N$ & $\displaystyle\max_{t\in [0,1]}\|e_h(t)\|_{L^1(\mathbb{R})}$ & $\displaystyle\max_{t\in [0,1]}\|e_h(t)\|_{L^2(\mathbb{R})}$ & $\displaystyle\max_{t\in [0,1]}\|e_h(t)\|_{L^\infty(\mathbb{R})}$ \\ \hline
    \multirow{2}{*}{$1.0$}   & $401$ & $1.95135(-05)$ & $1.85085(-05)$ & $3.28857(-05)$ \\
                             & $801$ & $1.93980(-05)$ & $1.83485(-05)$ & $3.18351(-05)$ \\ \hline
    \multirow{2}{*}{$0.1$}   & $401$ & $3.08028(-06)$ & $2.85654(-06)$ & $6.15974(-06)$ \\
                             & $801$ & $5.63505(-07)$ & $6.33056(-07)$ & $1.34676(-06)$ \\ \hline
    \multirow{2}{*}{$0.01$}  & $401$ & $2.78292(-05)$ & $7.10473(-05)$ & $4.45540(-04)$ \\
                             & $801$ & $3.71877(-06)$ & $9.20586(-06)$ & $6.08589(-05)$ \\ \hline
    \multirow{2}{*}{$0.001$} & $401$ & $3.81534(-04)$ & $3.23433(-03)$ & $4.24687(-02)$ \\
                             & $801$ & $5.88072(-05)$ & $5.14934(-04)$ & $9.48594(-03)$ \\ \hline
  \end{tabular}
  \label{tab:max_errors}
\end{table}

Table \ref{tab:max_errors} presents the maximum of the error norms of the finite element solution with $\delta_t = 10^{-3}$ against the analytic solution on the time interval $t\in [0, 1]$. We observe that for the moderate diffusion coefficients $\nu = 1.0, 0.1, 0.01$ errors in all computed norms are less than $10^{-4}$ by using a mesh of $801$ vertices. For the smaller diffusion coefficient $\nu = 0.001$ there exists a clear loss of accuracy, which indicates the necessity of applying more refined meshes.

The proposed numerical scheme applies less refined meshes as the solution spreads out. This can lead to a loss of accuracy for large times. Unfortunately, the computation of the analytic solution is too expensive to track the error norms for large time intervals. Alternatively, we next investigate the accuracy of the numerical solution for large times by studing its asymptotic behavior.

\subsection{Solution properties and its asymptotic behavior}\label{subsec:solution_properties}

Here, we show that the numerical finite element solution of the Burgers' equation still preserving some important properties of the analytic solution on the real line.

\begin{figure}[h!]
    \centering
    \begin{subfigure}[b]{0.45\textwidth}
    \centering
        \includegraphics[width=\textwidth]{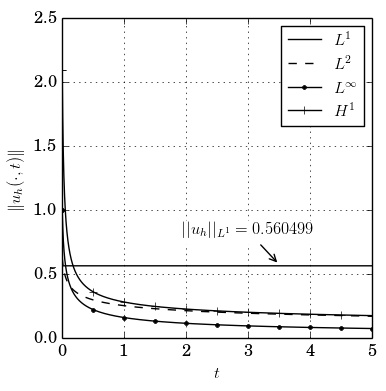}
        \caption{}
    \end{subfigure}
    ~
    \begin{subfigure}[b]{0.45\textwidth}
    \centering
        \includegraphics[width=\textwidth]{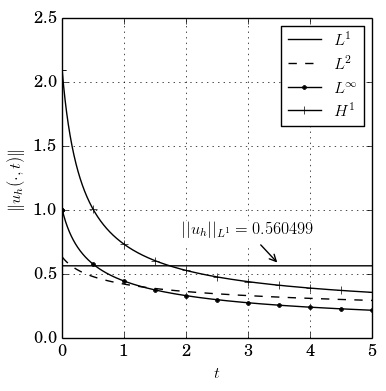}
        \caption{}
    \end{subfigure}
    \\
    \begin{subfigure}[b]{0.45\textwidth}
    \centering
        \includegraphics[width=\textwidth]{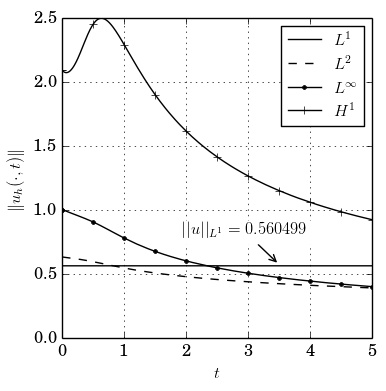}
        \caption{}
    \end{subfigure}
    ~
    \begin{subfigure}[b]{0.45\textwidth}
    \centering
        \includegraphics[width=\textwidth]{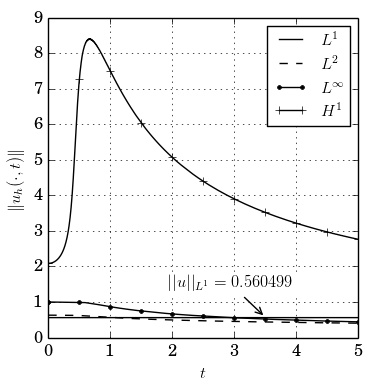}
        \caption{}
    \end{subfigure}
    \caption{Transient norms of the numerical results $u_h(\cdot,t)$ for the Burgers' equation with a mesh of $801$ vertices and diffusion coefficient: (a) $\nu = 1.0$, (b) $\nu = 0.1$, (c) $\nu = 0.01$, (d) $\nu = 0.001$.}\label{fig:norms-dt-1e-3-vert=801}
\end{figure}

Let's start by discussing on Figure \ref{fig:norms-dt-1e-3-vert=801}, which presents the profiles of the $L^1$, $L^2$, $L^\infty$ and $H^1$-norms of the numerical solution $u_h(\cdot,t)$ for $\nu = 1.0, 0.1, 0.01, 0.001$ with a mesh of $801$ vertices. The norms are computed as before in \eqref{eq:norms} and $H^1$-norm as follows:
\begin{equation}
  \|u_h\|_{H^1(\mathbb{R}} = \left[L \sum_{\mathcal{K}\in\mathcal{T}_h} \int_\mathcal{K} \left( |u_h(y,t)|^2 + \left|\frac{\p u_h(y,t)}{\p y}\right|^2\right)\,dy\right]^{\frac{1}{2}}.
\end{equation}

One of the most import characteristics of the Burgers' equation on the real line is that its solution conserves mass. The proposed second order finite element scheme is not conservative, nevertheless we observe that the numerical solutions have a small loss of mass, since the $L^1$-norm of the solution are constant up to five significant digits (see Figure \ref{fig:norms-dt-1e-3-vert=801}). This is in accordance with the accuracy reported in the last section.

Another observable characteristic of the numerical solution is that the $L_2$- and $L_\infty$-norms monotonically decreases as time increases. For $\nu = 1.0$ and $0.1$ the solution is highly diffusive causing the $H_1$-norm also to decrease with the time (see Figures \ref{fig:norms-dt-1e-3-vert=801}(a) and (b)).  But, for $\nu=0.01$ and $0.001$ convective effects is stronger than diffusive effects for small times, implying the increasing of the $H_1$-norm. However, as time increases so are the diffusion effects, and the $H_1$-norm become again a monotonically decreasing function (see Figures \ref{fig:norms-dt-1e-3-vert=801}(c) and (d)).

We finish this section by comparing the asymptotic behavior of the numerical against the analytic solutions. More specifically, we check if our numerical solution reproduces the solution behavior for large times described by equation \eqref{eq:asymptote}. In order to fix notation, we define:
\begin{equation}
 \tilde{\gamma}_{p} := t_\infty^{\frac{1}{2}\left(1 - \frac{1}{p}\right)}\|u_h(\cdot, t)\|_{L^p(\mathbb{R})},\quad p = 1, 2, \infty,
\end{equation}
where $t_\infty$ is such that $\|u_h(\cdot,t_f) - u_h(\cdot,t_f-\delta_t)\| < 10^{-10}$, where $\delta_t = 0.128$, and $\|\cdot\|$ denotes the $l_2-$vector norm. 

Table \ref{tab:asymptote} shows the comparison of the analytic $\gamma_p$ and numerical $\tilde{\gamma_p}$ for different values of $\nu$, mesh size, and $p=1, 2$ and $\infty$. We observe that the asymptotic of the numerical solutions agree with at least $3$ significant digits with the analytic solutions for all considered diffusion coefficients when a mesh of $801$ vertices is applied. These fine tests indicate that the proposed numerical scheme is accurate also for large times.

\begin{table}[h!]
  \centering
  \caption{Asymptotic behavior of the numerical \emph{versus} analytic solutions.}
  \begin{tabular}[h!]{lrccc} \hline
    $\nu$ & $N$ & $\gamma_p$ & $\tilde{\gamma}_p$ \\ \hline
    \multicolumn{4}{c}{$p=1$} \\ \hline
    \multirow{2}{*}{$1.0$}   & $401$  & \multirow{2}{*}{$5.60499(-01)$} & $5.60499(-01)$ \\
                             & $801$  & & $5.60499(-01)$ \\ \hline
    \multirow{2}{*}{$0.1$}   & $401$  & \multirow{2}{*}{$5.60499(-01)$} & $5.60499(-01)$ \\
                             & $801$  & & $5.60499(-01)$ \\ \hline
    \multirow{2}{*}{$0.01$}  & $401$  & \multirow{2}{*}{$5.60499(-01)$} & $5.60499(-01)$ \\
                             & $801$  & & $5.60499(-01)$ \\ \hline
    \multirow{2}{*}{$0.001$} & $401$  & \multirow{2}{*}{$5.60499(-01)$} & $9.25328(-01)$ \\
                             & $801$  & & $5.60509(-01)$ \\ \hline
    \multicolumn{4}{c}{$p=2$} \\ \hline
    \multirow{2}{*}{$1.0$}   & $401$  & \multirow{2}{*}{$2.50288(-01)$} & $2.50290(-01)$ \\
                             & $801$  & & $2.50292(-01)$ \\ \hline
    \multirow{2}{*}{$0.1$}   & $401$  & \multirow{2}{*}{$4.38152(-01)$} & $4.38153(-01)$ \\
                             & $801$  & & $4.38154(-01)$ \\ \hline
    \multirow{2}{*}{$0.01$}  & $401$  & \multirow{2}{*}{$5.92341(-01)$} & $5.92348(-01)$ \\
                             & $801$  & & $5.92351(-01)$ \\ \hline
    \multirow{2}{*}{$0.001$} & $401$  & \multirow{2}{*}{$6.23646(-01)$} & $6.23728(-01)$ \\
                             & $801$  & & $6.23691(-01)$ \\ \hline
    \multicolumn{4}{c}{$p=\infty$} \\ \hline
    \multirow{2}{*}{$1.0$}   & $401$  & \multirow{2}{*}{$1.58067(-01)$} & $1.58070(-01)$ \\
                             & $801$  & & $1.58071(-01)$ \\ \hline
    \multirow{2}{*}{$0.1$}   & $401$  & \multirow{2}{*}{$4.86580(-01)$} & $4.86575(-01)$ \\
                             & $801$  & & $4.86582(-01)$ \\ \hline
    \multirow{2}{*}{$0.01$}  & $401$  & \multirow{2}{*}{$9.25328(-01)$} & $9.25328(-01)$ \\
                             & $801$  & & $9.25359(-01)$ \\ \hline
    \multirow{2}{*}{$0.001$} & $401$  & \multirow{2}{*}{$1.03902(+00)$} & $1.04146(+00)$ \\
                             & $801$  & & $1.03968(+00)$  \\ \hline
  \end{tabular}
  \label{tab:asymptote}
\end{table}

\section{Final considerations}\label{sec:final_considerations}

In this paper we have presented a simple and accurate second order finite element scheme to simulate the Burgers' equation defined on the whole real line and subjected to initial conditions with compact support. The applied numerical scheme takes advantage of the convective-diffusive properties of this equation, which allow us to simulate the problem as a sequence of discrete auxiliary homogeneous Dirichlet's problems with fixed meshes.

Direct comparisons between analytic and numerical solutions have shown that the proposed scheme has good accuracy. Also, we have seen that the obtained numerical solutions reproduce that asymptotic behavior of the analytic solutions, which indicate that the simulations keep a good accuracy on very large time intervals. 

We close by observing that the good accuracy and the generality of the developed approach makes it suitable to produce insights about analytical properties of the solutions of the Burgers' equation. By appropriate modifications it may be suitable for related equations, for instance, for Burgers' equation in heterogeneous media.

\bibliographystyle{elsarticle-num}
\bibliography{main.bib}

\begin{thebibliography}{10}
\expandafter\ifx\csname url\endcsname\relax
  \def\url#1{\texttt{#1}}\fi
\expandafter\ifx\csname urlprefix\endcsname\relax\def\urlprefix{URL }\fi
\expandafter\ifx\csname href\endcsname\relax
  \def\href#1#2{#2} \def\path#1{#1}\fi

\bibitem{Bateman1915a}
M.~Bateman, Some recent researches on the motion of fluids, Mon. Wea. Rev. 43
  (1915) 163--170.

\bibitem{Burgers1974a}
J.~Burgers, The nonlinear diffusion equation, Springer, 1974.

\bibitem{Fletcher1982a}
F.~C.A., Numerical Solutions of Partial Differential Equations, Nort-Holland,
  Amsterdam, 1982, Ch. Burgers’ equation: a model for all reasons, pp.
  139--225.

\bibitem{Bastos2006a}
M.~Basto, V.~Semiao, F.~Calheiros, Dynamics in spectral solutions of burgers
  equation, J. Comput. Appl. Math. 205 (2006) 296--304.

\bibitem{Abd-el-Malek2000a}
M.~A. el~Malek, S.~El-Mansi, Group theoretic methods applied to burgers'
  equation, J. Comput. Appl. Math. 115 (2000) 1--12.

\bibitem{Evans2010a}
L.~Evans, Partial differential equations, 2nd Edition, Vol.~19 of Graduate
  Studies in Mathematics, The American Mathematical Society, 2010.

\bibitem{Gorguis2006a}
A.~Gorguis, A comparision between cole-hopf transformation and the
  decompisition method for solving burgers' equations, Appl. Math. Comput. 173
  (2006) 126--136.

\bibitem{Holland1977a}
C.~Holland, On the limiting behavior of burger's equation, J. Math. Anal. Appl.
  57 (1977) 156--160.

\bibitem{Rodin1970a}
E.~Rodin, On some approximate and exact solutions of boundary value problems
  for burgers' equation, J. Math. Anal. Appl. 30 (1970) 401--414.

\bibitem{Wood2006a}
W.~Wood, An exact solution for burger's equation, Commun. Numer. Meth. Engng.
  22 (2006) 797--798.

\bibitem{Aksan2005a}
E.~Aksan, A numerical solution of burgers' equation by finite element method
  constructed on the method of discretization in time, Appl. Math. Comput. 170
  (2005) 895--904.

\bibitem{Aksan2006a}
E.~Aksan, Quadratic {B}-spline finite element method for numerical solution of
  the burgers' equation, Appl. Math. Comput. 174 (2006) 884--896.

\bibitem{Arminjon1981a}
P.~Arminjon, C.~Beauchamp, Continuous and discontinuous finite element methods
  for burgers' equation, Comput. Methods Appl. Mech. Engrg. 25 (1981) 65--84.

\bibitem{Caldwell1982a}
J.~Caldwell, P.~Smith, Solution of burgers' equation with a large reynolds
  number, Appl. Math. Modelling 6 (1982) 381--385.

\bibitem{Caldwell1981a}
J.~Caldwell, P.~Wanless, A.~Cook, A finite element approach to burgers'
  equation, Appl. Math. Modelling 5 (1981) 189--193.

\bibitem{Caldwell1987a}
J.~Caldwell, P.~Wanless, A.~Cook, Solution of burgers' equation for large
  reynolds number using finite elements with moving nodes, Appl. Math.
  Modelling 11 (1987) 211--214.

\bibitem{Dogan2004a}
A.~Dogan, A galerkin finite element method to burgers' equation, Appl. Math.
  Comput. 157 (2004) 331--346.

\bibitem{Fletcher1983a}
C.~Fletcher, A comparison of finite element and finite difference solutions of
  the one- and two-dimensional burgers' equations, J. Comput. Phys. 51 (1983)
  159--188.

\bibitem{Hrymak1986a}
A.~Hrymak, G.~McRae, A.~Westerberg, An implementation of a moving finite
  element method, J. Comput. Phys. 63 (1986) 168--190.

\bibitem{Kadalbajoo2006a}
M.~Kadalbajoo, A.~Awasthi, A numerical method based on crank-nicolson scheme
  for burgers' equation, Appl. Math. Comput. 182 (2006) 1430--1442.

\bibitem{Kutluay2004a}
S.~Kutluay, A.~Esen, I.~Dag, Numerical solutions of the burgers' equation by
  the least-squares quadratic {B}-spline finite element method, J. Comput.
  Appl. Math. 167 (2004) 21--33.

\bibitem{Ozis2003a}
T.~\"{O}zi\c{s}, E.~Aksan, A.~\"{O}zde\c{s}, A finite element approach for
  solution of burgers' equation, Appl. Math. Comput. 139 (2003) 417--428.

\bibitem{Ozis2005b}
T.~\"{O}zi\c{s}, A.~Esen, S.~Kutluay, Numerical solution of burgers' equation
  by quadratic {B}-spline finite elements, Appl. Math. Comput. 165 (2005)
  237--249.

\bibitem{Shao2011a}
L.~Shao, X.~Feng, Y.~He, The local discontinuous {G}alerkin finite element
  method for burgers' equation, Math. Comput. Model. 54 (2011) 2943--2954.

\bibitem{Zhang2009a}
X.~Zhang, J.~Ouyang, L.~Zhang, Element-free characteristic galerkin method for
  burgers' equation, Eng. Anal. Boundary Elem. 33 (2009) 356--362.

\bibitem{Gulsu2006a}
M.~G\"{u}lsu, A finite difference approach for solution of burgers' equation,
  Appl. Math. Comput. 175 (2006) 1245--1255.

\bibitem{Kutluay1999a}
S.~Kutluay, A.~Bahadir, A.~\"{O}zde\c{s}, Numerical solution of one-dimensional
  burgers equation: explicit and exact-explicit finite difference methods, J.
  Comput. Appl. Math. 103 (1999) 251--261.

\bibitem{Mukundan2015a}
V.~Mukundan, A.~Awasthi, Efficient numerical techniques for burgers' equation,
  Appl. Math. Comput. 262 (2015) 282--297.

\bibitem{Aksan2004a}
E.~Aksan, A.~\"{O}zde{c}s, A numerical solution of burgers' equation, Appl.
  Math. Comput. 156 (2004) 395--402.

\bibitem{Caldwell1985a}
J.~Caldwell, R.~Saunders, P.~Wanless, A note on variation-iterative schemes
  applied to burgers' equation, J. Comput. Phys. 58 (1985) 275--281.

\bibitem{Ozis1996a}
T.~Ozis, A.~Ozdes, A direct variational methods applied to burgers' equation,
  J. Comput. Appl. Math. 71 (1996) 163--175.

\bibitem{Basdevant1986a}
C.~Basdevant, M.~Deville, P.~Haldenwang, J.~Lacroix, J.~Quazzani, R.~Peyret,
  P.~Orlandi, Spectral and finite difference solutions of the burgers'
  equation, Comput Fluids 14 (1986) 23--41.

\bibitem{Khater2008a}
A.~Khater, R.~Temsah, M.~Hassan, A {C}hebyshev spectral collocation method for
  solving burgers'-type equations, J. Comput. Appl. Math. 222 (2008) 333--350.

\bibitem{Hon1998a}
Y.~Hon, X.~Mao, An efficient numerical scheme for burgers' equation, J. Comput.
  Appl. Math. 95 (1998) 37--50.

\bibitem{Ozis2005a}
T.~\"{O}zi\c{s}, Y.~Aslan, The semi-approximate approach for solving burgers'
  equation with high reynolds number.

\bibitem{Tabatabaei2007a}
A.~Tabatabaei, E.~Shakour, M.~Dehghan, Some implicit methods for the numerical
  solution of burgers' equation, Appl. Math. Comput. 191 (2007) 560--570.

\bibitem{Inc2008a}
M.~Inc, On numerical solution of burgers' equation by homotopy analysis method,
  J. Phys. A 372 (2008) 356--360.

\bibitem{Saka2008a}
B.~Saka, I.~Da\u{g}, A numerical study of the burgers' equation, J. Frankl.
  Inst. 345 (2008) 328--348.

\bibitem{Hashemian2008a}
A.~Hashemian, H.~Shodja, A meshless approach for solution of burgers' equation,
  J. Comput. Appl. Math. 220 (2008) 226--239.

\bibitem{Xu2011a}
M.~Xu, R.-H. Wang, J.-H. Zhang, Q.~Fang, A novel numerical scheme for solving
  burgers' equation, Appl. Math. Comput. 217 (2011) 4473--4482.

\bibitem{Jiwari2015a}
R.~Jiwari, A hybrid numerical scheme for the numerical solution of the burgers'
  equation, Comput. Phys. Commun. 188 (2015) 59--67.

\bibitem{Cole1951a}
J.~D. Cole, et~al., On a quasi-linear parabolic equation occurring in
  aerodynamics, Quart. Appl. Math 9~(3) (1951) 225--236.

\bibitem{Hopf1950a}
E.~Hopf, The partial differential equation $u_t + uu_x = \mu u_{xx}$., Comm.
  Pure and Appl. Math. 3 (1950) 201--230.

\bibitem{Zingano1997a}
P.~Zingano, Some asymptotic limits for solutions of burgers equation, available
  at: \url{http://arxiv.org/pdf/math/0512503.pdf}, universidade Federal do Rio
  Grande do Sul (1997).

\bibitem{Johnson2009a}
C.~Johnson, Numerical solutions of partial differential equations by the finite
  element method, Dover, 2009.

\bibitem{dealII83}
W.~Bangerth, T.~Heister, L.~Heltai, G.~Kanschat, M.~Kronbichler, M.~Maier,
  B.~Turcksin, The \texttt{deal.II} library, version 8.3, preprint.

\bibitem{Scipy}
E.~Jones, T.~Oliphant, P.~Peterson, et~al., {Scipy} - open source scientific
  tool for {Python}: module for integration and {ODE}s, online, accessed on
  Mar/2016: \url{http://scipy.org/} (2001).

\end{thebibliography}

\end{document}